\documentclass{amsart}
\usepackage{amssymb}

\theoremstyle{plain}
\newtheorem{theorem}{Theorem}[section]

\newtheorem{lemma}[theorem]{Lemma}
\newtheorem{prop}[theorem]{Proposition}

\theoremstyle{remark}

\theoremstyle{definition}
\newtheorem{dfn}[theorem]{Definition}

\newtheorem*{example*}{Example}
\newtheorem*{examples*}{Examples}

\numberwithin{equation}{section}

\DeclareMathOperator{\Obj}{Obj}
\DeclareMathOperator{\Mor}{Mor}

\DeclareMathOperator{\lsp}{span}
\DeclareMathOperator{\supp}{supp}

\DeclareMathOperator{\Seg}{Seg}

\newcommand{\Lmin}{\Lambda^{\min}}

\newcommand{\field}[1]{\mathbb{#1}}
\newcommand{\CC}{\field{C}}

\newcommand{\NN}{\field{N}}

\newcommand{\TT}{\field{T}}
\newcommand{\ZZ}{\field{Z}}

\newcommand{\Gg}{{\mathcal G}}

\newcommand{\Kk}{{\mathcal K}}
\newcommand{\Ll}{{\mathcal L}}

\newcommand{\Oo}{{\mathcal O}}

\newcommand{\CE}{{\mathcal{CE}}}

\newcommand{\Tt}{{\mathcal T}}

\begin{document}

\title[topological higher-rank graphs]{Topological higher-rank
graphs and the $C^*$-algebras of topological $1$-graphs}%
\author{Trent Yeend}
\address{The University of Newcastle, Australia}
\email{trent.yeend@newcastle.edu.au}


\date{\today}

\begin{abstract}
We introduce the notion of a topological higher-rank graph, a
unified generalization of the higher-rank graph and the topological
graph. Using groupoid techniques, we define the Toeplitz and
Cuntz-Krieger algebras of topological higher-rank graphs, and show
that the $C^*$-algebras defined are coherent with the existing
theory.
\end{abstract}

\maketitle

\section{Introduction}\label{sec:intro}

We are interested in the generalization of directed-graph
$C^*$-algebras, of which, loosely speaking, there are two main
themes: either the graphs to which we associate $C^*$-algebras are
of higher-rank (for example, \cite{KP2,RSY1,RS,KP3,RSY2,S1,S2,FMY1,Ev}), or they
are given topological structure (for example, \cite{De1, KajWat1, De4, Ren4, DeKM, DeM, Kat1, MTom2, Kat2}). Our objective in this article is to present a unified approach to these two strands through the employment of groupoid theory.

This article explores a common approach to the $C^*$-algebras of higher-rank graphs and topological graphs. We begin by introducing topological higher-rank graphs -- a unified generalization of higher-rank graphs and topological graphs.

Given a topological higher-rank graph $\Lambda$, we define a path space $X_\Lambda$ and a topology on $X_\Lambda$. There is a natural action of $\Lambda$ on $X_\Lambda$, and from it we define a groupoid $G_\Lambda$ which has $X_\Lambda$ as its unit space; we call $G_\Lambda$ the path groupoid of $\Lambda$.

The topology on $X_\Lambda$ pulls back to give a topology on $G_\Lambda$. However, in general the topology may not be locally compact, and the groupoid's range and source maps may fail to be continuous. Hence we restrict our attention to the class of compactly aligned topological higher-rank graphs -- these are the topological analogues of the finitely aligned higher-rank graphs studied in \cite{RSY2,PatWel,FMY1}.

For a compactly aligned topological higher-rank graph $\Lambda$, $G_\Lambda$ is a locally compact topological groupoid which is $r$-discrete in the sense that the unit space $G_\Lambda^{(0)}$ is open in $G_\Lambda$. Furthermore, the range and source maps of $G_\Lambda$ are local homeomorphisms, so there is a Haar systems of counting measures on $G_\Lambda$. This allows us to define the full groupoid $C^*$-algebra $C^*(G_\Lambda)$, which we refer to as the Toeplitz algebra of $\Lambda$.

Identifying a closed invariant subset $\partial\Lambda$ of the unit space $G_\Lambda^{(0)}$, we define the boundary-path groupoid of $\Lambda$ to be $\Gg_\Lambda := G_\Lambda|_{\partial\Lambda}$ -- a locally compact $r$-discrete topological groupoid admitting a system of counting measures -- and define the Cuntz-Krieger algebra of $\Lambda$ to be $C^*(\Gg_\Lambda)$.

We finish by showing that when we restrict our attention to
topological graphs, we recover the same $C^*$-algebras as those
obtained in the existing theory.

Our work builds upon the theories of directed-graph $C^*$-algebras,
higher-rank graph $C^*$-algebras and topological-graph
$C^*$-algebras, so we mention some of the examples which these
classes contain.

The class of directed-graph $C^*$-algebras contains all nonunital,
simple, purely infinite, nuclear (SPIN) $C^*$-algebras whose $K_1$
group is torsion-free \cite{Sz3}, as well as examples of
$C^*$-algebras of quantum spaces \cite{HajMatSz, HonSz}. Up to
Morita equivalence, graph $C^*$-algebras provide all AF-algebras
\cite{Dr, Tyl} and unital SPIN $C^*$-algebras with torsion-free
$K_1$ \cite{Sz3}.

The class of higher-rank graph $C^*$-algebras provides more examples
of SPIN $C^*$-algebras, possibly with non-torsion-free $K_1$
\cite{RobSt, Ev}. Up to Morita equivalence, the class provides
examples of simple A$\TT$-algebras with real-rank 0 such as
irrational rotation algebras and Bunce Deddens algebras \cite{PRS};
in particular, these $C^*$-algebras are simple but neither AF nor
purely infinite.

Topological graphs generalize directed graphs as well as partially
defined local homeomorphisms on locally compact Hausdorff spaces.
Katsura shows in \cite{Kat2, Kat4} that the class of
topological-graph $C^*$-algebras contains all AF-algebras and many
AH-algebras, purely infinite $C^*$-algebras and stably
projectionless $C^*$-algebras. Katsura also shows that from
topological graphs one can obtain Exel-Laca algebras (see also
\cite{Schw}), ultragraph algebras, Matsumoto algebras, homeomorphism
$C^*$-algebras such as crossed products by partial homeomorphisms,
$C^*$-algebras associated with branched coverings, and
$C^*$-algebras associated with singly generated dynamical systems.

\subsection*{Higher-rank graphs}
Higher-rank graphs, or $k$-graphs ($k$ being an element of the
natural numbers and the given graph's rank), were first introduced
in \cite{KP2} as a unified way of approaching the higher-rank
Cuntz-Krieger algebras studied by Robertson and Steger \cite{RobSt}
and the Cuntz-Krieger algebras of directed graphs. The Toeplitz and
Cuntz-Krieger algebras of higher-rank graphs are collectively known
as higher-rank graph $C^*$-algebras.

Directed-graph $C^*$-algebras are recovered as the $C^*$-algebras of $1$-graphs. However, to pass from a directed graph $E$ to the
corresponding $1$-graph, we exchange the quadruple $E=(E^0,
E^1,r,s)$ for the pair $(E^*,l)$, where $E^*$ is the free category
generated by $E$, or the finite-path category of $E$, and
$l:E^*\to\NN$ is the functor which describes the length of each
path.

In general, a $k$-graph $(\Lambda,d)$ comprises a countable category
$\Lambda$, where morphisms are referred to as paths and objects as
vertices, together with a degree functor $d:\Lambda\to\NN^k$
describing the degree, or `shape', of each path.  The $k$-graph
$(\Lambda,d)$ is subject to one more condition: that of unique
factorizations of paths. Just as we may uniquely factorize a path in
a directed graph into its constituent edges, given a path
$\lambda\in\Lambda$ and any $m,n\in\NN^k$ such that $d(\lambda) =
m+n$, there are unique paths $\xi,\eta\in\Lambda$ such that $d(\xi)
= m$, $d(\eta) = n$ and $\lambda = \xi\eta$.

Given a $k$-graph $(\Lambda,d)$, one may represent it in a
$C^*$-algebra by a family of partial isometries $\{S_\lambda
:\lambda\in\Lambda\}$ satisfying conditions which encode the
structure of the graph and which give rise to a structure theory
analogous to that of the original namesake algebras (see
\cite[Definition~2.5]{RSY2}). Hence one may define and study the
Toeplitz and Cuntz-Krieger algebras of $(\Lambda,d)$.

There are a variety of approaches to the theory of higher-rank graph
$C^*$-algebras. In \cite{KP2}, Kumjian and Pask associate a groupoid
to each $k$-graph, and define the Cuntz-Krieger algebra of the graph
to be the corresponding groupoid $C^*$-algebra; this groupoid model
is an extension of the groupoid approach to the Cuntz-Krieger
algebra of a directed graph used in \cite{KPRRen,KPR}. Due to
groupoid considerations, the $k$-graphs studied in \cite{KP2} are
row-finite and have no sources.

In \cite{RSY1}, the authors generalize the work of Kumjian and Pask
to locally convex row-finite $k$-graphs (possibly with sources) by
using a direct $C^*$-algebraic analysis on the Cuntz-Krieger
representations of the $k$-graphs. The same authors further
generalize the theory to the setting of finitely aligned higher-rank
graphs in \cite{RSY2}, again using a direct analysis.

In \cite{FMY1} (see also \cite{PatWel}), the theory is brought
full-circle with an inverse semigroup and groupoid approach to the
Toeplitz and Cuntz-Krieger algebras of finitely aligned higher-rank
graphs; the approach recovers the Kumjian-Pask groupoid in the
`row-finite and no sources' setting, and is in many respects an
extension of the inverse semigroup and groupoid approach to
directed-graph $C^*$-algebras taken by Paterson in \cite{Pat2}.

Product systems of Hilbert bimodules may also be employed to
approach the Toeplitz algebras of finitely aligned higher-rank
graphs \cite{RS}, although at present there is a problem in the
extending of methods to the theory of Cuntz-Krieger algebras; this
intriguing problem is tied up with the need for an appropriate
notion of Cuntz-Pimsner covariance of product-system representations
for which there is a tractable structure theory.

\subsection*{Topological graphs}

A topological graph $E$ is a quadruple $(E^0,E^1,r,s)$, where $E^0$
and $E^1$ are locally compact Hausdorff spaces, $r,s:E^1\to E^0$ are
continuous, and $s$ is a local homeomorphism. Postponing details
until Section~\ref{sec:algebras of 1-graphs}, we give an outline of
a process by which one may associate Toeplitz and Cuntz-Krieger
algebras to each topological graph.

Given a topological graph $E$, Katsura \cite{Kat1} constructs a
right-Hilbert $C_0(E^0)-C_0(E^0)$ bimodule $C_d(E^1)$, where
$C_d(E^1)$ is a completion of $C_c(E^1)$. The Toeplitz algebra of
$E$ is then defined to be the universal $C^*$-algebra for Toeplitz,
or isometric, representations of $C_d(E^1)$, and the Cuntz-Krieger
algebra of $E$ is defined to be the $C^*$-algebra universal for
Cuntz-Pimsner covariant, or fully co-isometric, representations of
$C_d(E^1)$. The Toeplitz and Cuntz-Krieger algebras of topological
graphs are known collectively as topological-graph $C^*$-algebras.

Regarding directed graphs as topological graphs with discrete second
countable topologies, one recovers the same Hilbert bimodules and
$C^*$-algebras as those of the existing theory (see
\cite[Example~1.2]{FowR2} and \cite[Example~1 of Section~2]{Kat1}).

The direct $C^*$-algebraic analysis used with much success in the
discrete setting does not work for topological graphs: the direct
analysis involves manipulating generating families of partial
isometries and projections, whereas for a general topological graph
there may be no such partial isometries or projections present in
the representations.

Topological-graph $C^*$-algebras may also be approached using
groupoid methods \cite{ArRen,Ren4}, \cite[Section~10.3]{Kat2}. The
groupoid approach was first achieved by Deaconu \cite{De1} for
topological graphs with compact vertex and edge sets, with $r$ a
homeomorphism, and with $s$ a surjective local homeomorphisms.

\section{Topological higher-rank graphs}\label{sec:top h-r graphs}

Given $k\in\NN$, a \emph{topological $k$-graph} is a pair
$(\Lambda,d)$ consisting of a small category
$\Lambda=(\Obj(\Lambda),\Mor(\Lambda),r,s,\circ)$ and a functor
$d:\Lambda\to\NN^k$, called the \emph{degree map}, which satisfy:
\begin{enumerate}
\item
$\Obj(\Lambda)$ and $\Mor(\Lambda)$ are second-countable locally
compact Hausdorff spaces;
\item
$r,s:\Mor(\Lambda)\to\Obj(\Lambda)$ are continuous and $s$ is a
local homeomorphism;
\item
Composition $\circ : \Lambda\times_c\Lambda\to\Lambda$ is continuous
and open;
\item
$d$ is continuous, where $\NN^k$ has the discrete topology;
\item
For all $\lambda\in\Lambda$ and $m,n\in\NN^k$ such that
$d(\lambda)=m+n$, there exists unique
$(\xi,\eta)\in\Lambda\times_c\Lambda$ such that $\lambda=\xi\eta$,
$d(\xi)=m$ and $d(\eta)=n$.
\end{enumerate}
We refer to the morphisms of $\Lambda$ as \emph{paths} and to the
objects of $\Lambda$ as \emph{vertices}. The codomain and domain
maps in $\Lambda$ are called the range and source maps,
respectively.

We use the partial ordering of $\NN^k$,
\[
m\le n \iff m_i\le n_i \text{ for } i=1,\dots,k,
\]
and use the notation $\vee$ and $\wedge$ for the coordinate-wise
maximum and minimum.

For $m\in\NN^k$, define $\Lambda^m$ to be the set $d^{-1}(\{m\})$ of
paths of degree $m$. Define $\Lambda *_s\Lambda
:=\{(\lambda,\mu)\in\Lambda\times\Lambda : s(\lambda)=s(\mu)\}$, and
for $U,V\subset\Lambda$ define $U*_sV:=(U\times V)\cap (\Lambda
*_s\Lambda)$ and $UV:= \{\lambda\mu : (\lambda,\mu)\in U\times_c
V\}$; in particular, for $v\in\Lambda^0$, $vU :=\{v\}U= \{\lambda\in
U : r(\lambda)=v\}$ and similarly $Uv := \{\lambda\in U :
s(\lambda)=v\}$. For $p,q\in\NN^k$, $U\subset\Lambda^p$ and
$V\subset\Lambda^q$, we write
\[
U\vee V := U\Lambda^{(p\vee q)-p}\cap V\Lambda^{(p\vee q)-q}
\]
for the set of \emph{minimal common extensions} of paths from $U$
and $V$. For $\lambda,\mu\in\Lambda$, we write
\[
\Lmin(\lambda,\mu):=\{(\alpha,\beta) : \lambda\alpha=\mu\beta,
d(\lambda\alpha)=d(\lambda)\vee d(\mu)\}
\]
for the set of pairs which give minimal common extensions of
$\lambda$ and $\mu$; that is,
\[
\Lmin(\lambda,\mu) = \{(\alpha,\beta) : \lambda\alpha=\mu\beta\in
\{\lambda\}\vee\{\mu\}\}.
\]

\section{The path groupoid}

To describe our path groupoid $G_\Lambda$, we need some terminology.
Let $(\Lambda_1,d_1)$ and $(\Lambda_2,d_2)$ be topological
$k$-graphs.  A \emph{graph morphism} between $\Lambda_1$ and
$\Lambda_2$ is a continuous functor $x:\Lambda_1\to\Lambda_2$
satisfying $d_2(x(\lambda))=d_1(\lambda)$ for all
$\lambda\in\Lambda_1$.

For $k\in\NN$ and $m\in (\NN\cup\{\infty\})^k$, define the
topological $k$-graph $(\Omega_{k,m},d)$ by giving the discrete
topologies to the sets
\[
\Obj(\Omega_{k,m}):=\{p\in\NN^k : p\le m\}
\]
and
\[
\Mor(\Omega_{k,m}):= \{(p,q)\in\NN^k\times\NN^k : p\le q\le m\},
\]
and setting $r(p,q):=p$, $s(p,q):=q$, $(n,p)\circ (p,q):= (n,q)$ and
$d(p,q):= q-p$.

Let $(\Lambda,d)$ be a topological $k$-graph. We define the
\emph{path space of $\Lambda$} to be
\[
X_\Lambda
:=\bigcup_{m\in(\NN\cup\{\infty\})^k}\{x:\Omega_{k,m}\to\Lambda : x
\text{ is a graph morphism}\}.
\]
We extend the range and degree maps to $x:\Omega_{k,m}\to\Lambda$ in
$X_\Lambda$ by setting $r(x):=x(0)$ and $d(x):=m$. For
$v\in\Lambda^0$ we define $v X_\Lambda:=\{x\in X_\Lambda :
r(x)=v\}$.

For $x\in X_\Lambda$, $m\in\NN^k$ with $m\le d(x)$, and
$\lambda\in\Lambda$ with $s(\lambda)=r(x)$, there exist unique graph morphisms $\lambda x$ and $\sigma^m x$ in $X_\Lambda$ satisfying $d(\lambda x)=d(\lambda)+d(x)$, $d(\sigma^m x) = d(x)-m$,
\[
(\lambda x)(0,p) = \begin{cases}
\lambda(0,p) &\text{if } p\le d(\lambda) \\
\lambda x(0,p-d(\lambda)) &\text{if } d(\lambda)\le p\le d(\lambda
x),
\end{cases}
\]
and
\[
(\sigma^m x)(0,p) = x(m,m+p) \quad\text{for } p\le d(\sigma^m x).
\]

For each $\lambda\in\Lambda$ there is a unique graph morphism
$x_\lambda:\Omega_{k,d(\lambda)}\to\Lambda$ such that
$x_\lambda(0,d(\lambda))=\lambda$; in this sense, we may view
$\Lambda$ as a subset of $X_\Lambda$ and we refer to elements of
$X_\Lambda$ as paths. Indeed, for $\lambda\in\Lambda$ and
$p,q\in\NN^k$ with $0\le p\le q\le d(\lambda)$ we write
$\lambda(0,p)$, $\lambda(p,q)$ and $\lambda(q,d(\lambda))$ for the
unique elements of $\Lambda$ which satisfy
$\lambda=\lambda(0,p)\lambda(p,q)\lambda(q,d(\lambda))$,
$d(\lambda(0,p))=p$, $d(\lambda(p,q))=q-p$ and
$d(\lambda(q,d(\lambda)))=d(\lambda)-q$.

\begin{dfn}
Let $(\Lambda,d)$ be a topological $k$-graph. The \emph{path
groupoid} $G_\Lambda$ has object set $\Obj(G_\Lambda):= X_\Lambda$,
morphism set
\begin{align*}
\Mor(G_\Lambda) &:= \{ (\lambda x, d(\lambda)-d(\mu), \mu x) \in
X_\Lambda\times\ZZ^k\times X_\Lambda: \\
&\phantom{MMMMMMNMMN}(\lambda,\mu)\in\Lambda *_s\Lambda,~ x\in
X_\Lambda \text{ and } s(\lambda)=r(x)\} \\
&=\{ (x,m,y) \in X_\Lambda\times\ZZ^k\times X_\Lambda : \text{ there
exist } p,q\in\NN^k \text{ such that}\\
&\phantom{MMMNNMMN} p\le d(x),~ q\le d(y),~ p-q=m \text{ and }
\sigma^p x = \sigma^q y \},
\end{align*}
range and source maps $r(x,m,y):= x$ and $s(x,m,y):= y$,
composition
\[
((x,m,y),(y,n,z))\mapsto (x,m+n,z),
\]
and inversion $(x,m,y)\mapsto (y,-m,x)$.
\end{dfn}

To define a topology on $G_\Lambda$, we use the following notation.
For $F\subset \Lambda *_s\Lambda$ and $m\in\ZZ^k$, define
$Z(F,m)\subset G_\Lambda$ by
\[
Z(F,m) := \{(\lambda x,d(\lambda)-d(\mu),\mu x)\in G_\Lambda :
(\lambda,\mu)\in F, d(\lambda)-d(\mu)=m\}.
\]
For $U\subset\Lambda$, define $Z(U)\subset G_\Lambda^{(0)}$ by
\[
Z(U) := Z(U*_sU,0)\cap Z(\Lambda^0*_s\Lambda^0,0).
\]

\begin{prop}[{\cite[Proposition~3.8]{Y1}}]\label{prop:defining a
topology on G_Lambda}
Let $(\Lambda,d)$ be a topological $k$-graph. The family of sets of
the form
\[
Z(U*_s V,m)\cap Z(F,m)^c,
\]
where $m\in\ZZ^k$, $U,V\subset\Lambda$ are open and $F\subset
\Lambda *_s\Lambda$ is compact, is a basis for a second-countable
Hausdorff topology on $G_\Lambda$.
\end{prop}

Sounds perfect! Not quite, there's a problem: If $\Lambda$ is not
\emph{compactly aligned} in the sense that, for compact
$U\subset\Lambda^p$ and $V\subset\Lambda^q$, the set $U\vee V$ is
compact, then $G_\Lambda$ is neither locally compact nor a
topological groupoid -- the range and source maps in $G_\Lambda$
fail to be continuous. The property of being compactly aligned is to
topological higher-rank graphs what finitely aligned is to discrete
higher-rank graphs (see \cite{RSY2,PatWel,FMY1}), and regardless of
the approach to the $C^*$-algebras of higher-rank graphs, the
property has presented itself as a necessary assumption. Also, it is
straightforward to see that every topological $1$-graph is compactly
aligned.

\begin{theorem}[{\cite[Theorem~3.16]{Y1}}]
Let $(\Lambda,d)$ be a compactly aligned topological $k$-graph. Then
$G_\Lambda$ is a locally compact $r$-discrete topological groupoid
admitting a Haar system consisting of counting measures.
\end{theorem}

We define the \emph{Toeplitz algebra of $\Lambda$} to be the full
groupoid $C^*$-algebra $C^*(G_\Lambda)$.

For a discrete directed graph $E$, Paterson \cite{Pat2} defines an
inverse semigroup $S_E$ and an action of $S_E$ on $X_{E^*}$. The
topological groupoid $H_E$ is then defined as the groupoid of germs
of the action. Comparison of the topological and groupoid structures
reveals that $H_E$ and $G_{E^*}$ are isomorphic as topological
groupoids.

Given a discrete finitely aligned higher-rank graph $(\Lambda,d)$,
the authors of \cite{FMY1} define an $r$-discrete groupoid
$G_\Lambda^{\rm [FMY]}$ and show that its $C^*$-algebra is
isomorphic to the Toeplitz algebra of $\Lambda$
\cite[Theorem~5.9]{FMY1}. The two groupoids $G_\Lambda^{\rm [FMY]}$
and $G_\Lambda$ are isomorphic, so we see that our Toeplitz algebra
$C^*(G_\Lambda)$ is coherent with the existing theory for
higher-rank graphs.

\section{The boundary-path groupoid}

Let $(\Lambda,d)$ be a topological $k$-graph and let
$V\subset\Lambda^0$. A set $E\subset V\Lambda$ is \emph{exhaustive}
for $V$ if for all $\lambda\in V\Lambda$ there exists $\mu\in E$
such that $\Lmin(\lambda,\mu)\neq\emptyset$. For $v\in\Lambda^0$,
let $v\CE(\Lambda)$ denote the set of all compact sets
$E\subset\Lambda$ such that $r(E)$ is a neighbourhood of $v$ and $E$
is exhaustive for $r(E)$.

A path $x\in X_\Lambda$ is called a \emph{boundary path} if for all
$m\in\NN^k$ with $m\le d(x)$, and for all $E\in x(m)\CE(\Lambda)$,
there exists $\lambda\in E$ such that $x(m,m+d(\lambda))=\lambda$.
We write $\partial\Lambda$ for the set of all boundary paths in
$X_\Lambda$. For $v\in\Lambda^0$ and $V\subset\Lambda^0$, we define
$v\partial\Lambda = \{x\in\partial\Lambda : r(x) = v\}$ and
$V(\partial\Lambda) = \{x\in\partial\Lambda : r(x)\in V\}$.

The set of boundary paths of $\Lambda$ is a nonempty closed
invariant subset of $G_\Lambda^{(0)}$, so we can make the following
definition.

\begin{dfn}
Let $(\Lambda,d)$ be a compactly aligned topological $k$-graph. The
\emph{boundary-path groupoid} $\Gg_\Lambda$ is the reduction
$\Gg_\Lambda := G_\Lambda|_{\partial\Lambda}$, which is a locally
compact $r$-discrete topological groupoid admitting a Haar system
consisting of counting measures. The \emph{Cuntz-Krieger algebra of
$\Lambda$} is the full groupoid $C^*$-algebra $C^*(\Gg_\Lambda)$.
\end{dfn}

For a discrete finitely aligned higher-rank graph $(\Lambda,d)$,
straightforward comparisons with the work in \cite{FMY1} show that
the Cuntz-Krieger algebra $C^*(\Lambda)$ is isomorphic to
$C^*(\Gg_\Lambda)$.

\section{$C^*$-algebras of topological
$1$-graphs}\label{sec:algebras of 1-graphs}

Let $E$ be a second-countable topological graph as defined in
\cite[Definition~2.1]{Kat1}; that is, $E = (E^0,E^1,r,s)$ is a
directed graph with $E^0$, $E^1$ second-countable locally compact
Hausdorff spaces, $r,s:E^1\to E^0$ continuous, and $s$ a local
homeomorphism. The free category generated by $E$, endowed with the
relative topology inherited from the union of the product
topologies, together with the length functor $l(e_1\cdots e_n):= n$,
forms a topological $1$-graph $(E^*,l)$. Conversely, given a
topological $1$-graph $(\Lambda,d)$, the quadruple $E_\Lambda:=
(\Lambda^0,\Lambda^1,r|_{\Lambda^1},s|_{\Lambda^1})$ is a
second-countable topological graph with $((E_\Lambda)^*,l)\cong
(\Lambda,d)$.

In this section, we will see that the $C^*$-algebras,
$C^*(G_\Lambda)$ and $C^*(\Gg_\Lambda)$, are isomorphic to the
Toeplitz and Cuntz-Krieger algebras, $\Tt(E_\Lambda)$ and
$\Oo(E_\Lambda)$, of the associated topological graph, as defined in
\cite{Kat1}.

We note that Muhly and Tomforde \cite{MTom2} elegantly generalize
the theory set forth by Katsura, removing the hypothesis that the
source map of the topological graph is a local homeomorphism; in its stead, the authors impose the weaker condition that the source map
is open and that there is a family of Radon measures
$\{\lambda_v\}_{v\in E^0}$ on $E^1$ satisfying:
\begin{enumerate}
\item
$\supp(\lambda_v)=s^{-1}(v)$ for all $v\in E^0$
\item
$v\mapsto \int_{E^1}\xi(e)d\lambda_v(e)$ is in $C_c(E^0)$ for all
$\xi\in C_c(E^1)$.
\end{enumerate}

\begin{theorem}\label{thm:Toeplitz algebra of top graphs
correspond}
For a topological $1$-graph, we have $\Tt(E_\Lambda)\cong
C^*(G_\Lambda)$.
\end{theorem}

\begin{theorem}\label{thm:C-K algebras of top graphs correspond}
For a topological $1$-graph, we have $\Oo(E_\Lambda)\cong
C^*(\Gg_\Lambda)$.
\end{theorem}

Henceforth $(\Lambda,d)$ will be a topological $1$-graph with
associated topological graph $E_\Lambda=(\Lambda^0,\Lambda^1,r,s)$.
We begin by defining the Hilbert bimodule
$_{C_0(\Lambda^0)}C_d(\Lambda^1)_{C_0(\Lambda^0)}$.

For $\xi\in C(\Lambda^1)$, define
$\langle\xi,\xi\rangle:\Lambda^0\to [0,\infty]$ by
\[
\langle\xi,\xi\rangle(v) := \sum_{e\in\Lambda^1v}|\xi(e)|^2,
\]
and define
\[
C_d(\Lambda^1) := \{\xi\in C(\Lambda^1) : \langle\xi,\xi\rangle\in
C_0(\Lambda^0)\}.
\]
For $\xi,\eta\in C_d(\Lambda^1)$, define $\langle\xi,\eta\rangle\in
C_0(\Lambda^0)$ by
\[
\langle\xi,\eta\rangle(v) :=
\sum_{e\in\Lambda^1v}\overline{\xi(e)}\eta(e),
\]
and define left and right actions of $C_0(\Lambda^0)$ on
$C_d(\Lambda^1)$ by
\[
(f\cdot \xi)(e) = (\phi(f)\xi)(e) := f(r(e))\xi(e)
\quad\text{and}\quad (\xi\cdot f)(e) := \xi(e)f(s(e)).
\]
Then $C_d(\Lambda^1)$ is a right-Hilbert $C_0(\Lambda^0)$-module and
$\phi$ is a homomorphism of $C_0(\Lambda^0)$ into
$\Ll(C_d(\Lambda^1))$. Hence $C_d(\Lambda^1)$ is a right-Hilbert
$C_0(\Lambda^0) - C_0(\Lambda^0)$-bimodule
\cite[Proposition~1.10]{Kat1}, and $C_c(\Lambda^1)$ is dense in
$C_d(\Lambda^1)$ \cite[Lemma~1.6]{Kat1}.

A \emph{Toeplitz $E_\Lambda$-pair} on a $C^*$-algebra $B$ is a pair
of maps $\Psi = (\Psi_0,\Psi_1)$ such that $\Psi_0:
C_0(\Lambda^0)\to B$ is a homomorphism and $\Psi_1:C_d(\Lambda^1)\to
B$ is a linear map satisfying
\begin{itemize}
\item[(i)]
$\Psi_1(\xi)^*\Psi_1(\eta) = \Psi_0(\langle\xi,\eta\rangle)$ for all
$\xi,\eta\in C_d(\Lambda^1)$, and
\item[(ii)]
$\Psi_0(f)\Psi_1(\xi) = \Psi_1(\phi(f)\xi)$ for all $f\in
C_0(\Lambda^0)$ and $\xi\in C_d(\Lambda^1)$.
\end{itemize}
The \emph{Toeplitz algebra of $E_\Lambda$}, denoted
$\Tt(E_\Lambda)$, is the $C^*$-algebra universal for Toeplitz
$E_\Lambda$-pairs.

It is straightforward to see that $\Psi_1$ is continuous and that
$\Psi_1(\xi)\Psi_0(f) = \Psi_1(\xi\cdot f)$ for all $\xi\in
C_d(\Lambda^1)$ and $f\in C_0(\Lambda^0)$.

Next we define three open subsets of $\Lambda^0$:
\[
\Lambda^0_{\rm sce} := \Lambda^0\setminus \overline{r(\Lambda^1)},
\]
\[
\Lambda^0_{\rm fin} := \{v\in\Lambda^0 : ~v\text{ has a
neighbourhood } V \text{ such that } V\Lambda^1 \text{ is
compact}\}
\]
and
\[
\Lambda^0_{\rm rg} := \Lambda^0_{\rm fin}\setminus
\overline{\Lambda^0_{\rm sce}}.
\]

Given a Toeplitz $E_\Lambda$-pair $\Psi$ on $B$, there is a
homomorphism $\Psi^{(1)}$ of $\Kk(C_d(\Lambda^1))$ into $B$
satisfying
\[
\Psi^{(1)}(\xi\otimes\eta^*) = \Psi_1(\xi)\Psi_1(\eta)^*,
\]
where $\xi\otimes\eta^*\in\Kk(C_d(\Lambda^1))$ is the `rank-one'
operator defined by $(\xi\otimes\eta^*)(\zeta) = \xi\cdot
\langle\eta,\zeta\rangle$ for all $\zeta\in C_d(\Lambda^1)$.

By \cite[Proposition~1.24]{Kat1}, the restriction of $\phi$ to
$C_0(\Lambda^0_{\rm rg})$ is an injection into
$\Kk(C_d(\Lambda^1))$, so we can make the following definition.

A \emph{Cuntz-Krieger $E_\Lambda$-pair} $\Psi = (\Psi_0,\Psi_1)$ is
a Toeplitz $E_\Lambda$-pair which satisfies:
\[
\Psi_0(f) = \Psi^{(1)}(\phi(f)) \quad\text{for all } f\in
C_0(\Lambda^0_{\rm rg}).
\]
The Cuntz-Krieger algebra of $E_\Lambda$, denoted $\Oo(E_\Lambda)$,
is the $C^*$-algebra universal for Cuntz-Krieger $E_\Lambda$-pairs.

We now construct a Toeplitz $E_\Lambda$-pair $\Psi$ on
$C^*(G_\Lambda)$. For $f\in C_0(\Lambda^0)$, define
$\psi_0(f):G_\Lambda^{(0)}\to\CC$ by
\[
\psi_0(f)(x,0,x) := f(r(x))=f(x(0)).
\]
For $f\in C_0(\Lambda^0)$, we have $\psi_0(f)\in
C_0(G_\Lambda^{(0)})$, and $\psi_0$ is a homomorphism.

Since $G_\Lambda$ is $r$-discrete, it follows that
$C_0(G_\Lambda^{(0)})$ is a subalgebra of $C^*(G_\Lambda)$. Letting
$\iota :C_0(G_\Lambda^{(0)})\to C^*(G_\Lambda)$ be the inclusion
homomorphism, we define $\Psi_0:C_0(\Lambda^0)\to C^*(G_\Lambda)$ by
$\Psi_0 = \iota\circ \psi_0$.

For $\xi\in C_c(\Lambda^1)$, define $\Psi_1(\xi):G_\Lambda\to\CC$
by
\[
\Psi_1(\xi)(x,m,y) = \delta_{m,1}\delta_{\sigma^1x,y}\xi(x(0,1)).
\]
Then $\Psi_1$ is linear, and for $\xi\in C_c(\Lambda^1)$, we have
$\Psi_1(\xi)\in C_c(G_\Lambda)$ with
\begin{equation}\label{eqn:support of Psi_1}
\supp(\Psi_1(\xi))\subset Z((\supp(\xi))*_s(s(\supp(\xi))),1).
\end{equation}

\begin{lemma}\label{lem:Toepltz relations on C_c}
For $f\in C_c(\Lambda^0)$ and $\xi,\eta\in C_c(\Lambda^1)$, we have
\begin{itemize}
\item[(i)]
$\Psi_1(\xi)^*\Psi_1(\eta) = \Psi_0(\langle\xi,\eta\rangle)$ and
\item[(ii)]
$\Psi_0(f)\Psi_1(\xi) = \Psi_1(\phi(f)\xi)$.
\end{itemize}
\end{lemma}

\begin{proof}
Since $\Psi_1(\xi)$, $\Psi_1(\eta)$ and $\Psi_0(f)$ have compact
support, we can use convolution on $C_c(G_\Lambda)$ to calculate
products. To begin,
\[
\Psi_1(\xi)^*\Psi_1(\eta)(x,m,y) = \sum_{(x,n,z)\in
G_\Lambda}\overline{\Psi_1(\xi)(z,-n,x)}\Psi_1(z,m-n,y).
\]
A summand on the right-hand side may be nonzero only if $-n=1$,
$\sigma^1z=x$, $m-n=1$ and $\sigma^1z=y$, which is precisely when
$m=0$, $x=y$, $n=-1$ and $z=ex$ for some $e\in\Lambda^1(r(x))$. So
the support of $\Psi_1(\xi)^*\Psi_1(\eta)$ is contained in
$G_\Lambda^{(0)}$, and for $(x,0,x)\in G_\Lambda^{(0)}$, we have
\begin{align*}
\Psi_1(\xi)^*\Psi_1(\eta)(x,0,x) &=
\sum_{e\in\Lambda^1r(x)}\overline{\Psi_1(\xi)(ex,1,x)}\Psi_1(ex,1,x)
\\
&= \sum_{e\in\Lambda^1r(x)}\overline{\xi(e)}\eta(e) \\
&= \langle \xi,\eta\rangle (r(x)) \\
&= \Psi_0(\langle\xi,\eta\rangle)(x,0,x),
\end{align*}
giving (i).

For (ii), we have
\[
\Psi_0(f)\Psi_1(\xi)(x,m,y) = \sum_{(x,n,z)\in
G_\Lambda}\Psi_0(f)(x,n,z)\Psi_1(\xi)(z,m-n,y).
\]
A summand on the right-hand side may be nonzero only if $n=0$,
$x=z$, $m-n=1$ and $\sigma^1z=y$. Hence the sum reduces to a single
summand, and we have
\begin{align*}
\Psi_0(f)\Psi_1(\xi)(x,m,y) &=
\delta_{m,1}\delta_{\sigma^1x,y}\Psi_0(f)(x,0,x)\Psi_1(x,1,\sigma^1x)
\\
&= \delta_{m,1}\delta_{\sigma^1x,y}f(r(x))\xi(x(0,1)) \\
&=\delta_{m,1}\delta_{\sigma^1x,y}(\phi(f)\xi)(x(0,1)) \\
&= \Psi_1(\phi(f)\xi)(x,m,y),
\end{align*}
as required.
\end{proof}

Using Lemma~\ref{lem:Toepltz relations on C_c}(i), we see that
$\Psi_1$ is norm-decreasing, and since $C_c(\Lambda^1)$ is dense in
$C_d(\Lambda^1)$, it follows that $\Psi_1$ extends to a
norm-decreasing linear map $\Psi_1:C_d(\Lambda^1)\to
C^*(G_\Lambda)$. Continuity allows us to extend the
properties in Lemma~\ref{lem:Toepltz relations on C_c} to all of
$C_0(\Lambda^0)$ and $C_d(\Lambda^1)$, giving:

\begin{prop}\label{prop:Toeplitz relations}
For $f\in C_0(\Lambda^0)$ and $\xi,\eta\in C_d(\Lambda^1)$,
\begin{itemize}
\item[(i)]
$\Psi_1(\xi)^*\Psi_1(\eta) = \Psi_0(\langle\xi,\eta\rangle)$ and
\item[(ii)]
$\Psi_0(f)\Psi_1(\xi) = \Psi_1(\phi(f)\xi)$.
\end{itemize}
\end{prop}

Therefore $(\Psi_0,\Psi_1)$ is a Toeplitz $E_\Lambda$-pair on
$C^*(G_\Lambda)$, so the universal property of $\Tt(E_\Lambda)$
gives a homomorphism $\Psi_0\times_\Tt\Psi_1:\Tt(E_\Lambda)\to
C^*(G_\Lambda)$.

The following notation is handy. For $m,p,q\in\NN$ with $p\le q\le
m$, define the continuous map
$\Seg^m_{(p,q)}:\Lambda^m\to\Lambda^{q-p}$ by
$\Seg^m_{(p,q)}(\lambda):=\lambda(p,q)$.

\begin{prop}\label{prop:surjectivity}
$\Psi_0\times_\Tt\Psi_1:\Tt(E_\Lambda)\to C^*(G_\Lambda)$ is
surjective.
\end{prop}

\begin{proof}
We complete the proof in two steps, first showing that
$C_0(G_\Lambda^{(0)})$ is in the image of $\Psi_0\times_\Tt\Psi_1$,
then using this to show $\Psi_0\times_\Tt\Psi_1$ is surjective.

Let $\{U_j\}_{j\in\NN}$ be a basis for $\Lambda^1$ comprising
relatively compact open sets $U_j$ such that $s|_{U_j}$ is a
homeomorphism. Define
\[
W_1 = \{\Psi_0(f) : f\in C_c(\Lambda^0)\},
\]
\[
\begin{split}
W_2 = \{\Psi_1(\xi)\cdots
\Psi_1(\xi_p)&\Psi_1(\eta_p)^*\cdots\Psi_1(\eta_1)^* : p\in\NN, \\
&\text{and each } \xi_i,\eta_i\in C_0(U_{j_i}) \text{ for some }
j_i\in\NN\}
\end{split}
\]
and
\[
W = \lsp (W_1\cup W_2).
\]
Then $W$ is a $*$-subalgebra of $C_0(G_\Lambda^{(0)})$, and $W$
separates points in $G_\Lambda^{(0)}$ and does not vanish
identically at any point of $G_\Lambda^{(0)}$. Therefore, by the
Stone-Weierstrass Theorem, $W$ is uniformly dense in
$C_0(G_\Lambda^{(0)})$. Since the supremum norm coincides with the
$I$-norm on elements of $C_0(G_\Lambda^{(0)})$, and the $I$-norm
bounds the $C^*$-norm, it follows that $W$ is dense in
$C_0(G_\Lambda^{(0)})$ with respect to the $C^*$-norm, and
$C_0(G_\Lambda^{(0)})$ is in the image of $\Psi_0\times_\Tt\Psi_1$.

Now fix $f\in C_c(G_\Lambda)$. Let $\{W_j\}_{j=1}^n$ be an open
cover of $\supp f$ of the form
\[
W_j = Z((U_1^j\cdots U_{p_j}^j)*_s(V_1^j\cdots V_{q_j}^j),p_j-q_j),
\]
where each $U_i^j$ and $V_{i'}^j$ are relatively compact open
subsets of $\Lambda^1$ such that $s|_{U_i^j}$ and $s|_{V_{i'}^j}$
are homeomorphisms. Let $\{\varphi_j\}_{j=1}^n$ be a partition of
unity subordinate to $\{W_j\}_{j=1}^n$; that is, each $\varphi_j :
G_\Lambda\to [0,1]$ is continuous and $\supp(\varphi_j)\subset W_j$,
and $\sum_{j=1}^n \varphi_j(x,m,y) =1$ for all $(x,m,y)\in
\supp(f)$. Then $f = \sum_j\varphi_jf$, where the products
$\varphi_jf$ are pointwise and $\supp(\varphi_jf)\subset W_j$.

Fix $j\in\{1,\dots,n\}$; we will show that $\varphi_jf$ is in the
image of $\Psi_0\times_\Tt\Psi_1$. For convenience, we drop the
index $j$, and write
\[
W_j = Z((U_1\cdots U_p)*_s(V_1\cdots V_q),p-q).
\]

For $x\in G_\Lambda^{(0)}$, it follows from injectivity of
$s|_{U_1\cdots U_p}$ and $s|_{V_1\cdots V_q}$ that there is at most
one $(\lambda,\mu)\in (U_1\cdots U_p)*_s(V_1\cdots V_q)$ such that
$s(\lambda) = r(x)$. Hence we may define a function
$g:G_\Lambda^{(0)}\to\CC$ by
\[
g(x) := \begin{cases}
\varphi_jf(\lambda x,p-q,\mu x) &\text{if there exists }
(\lambda,\mu) \text{ such that } (\lambda x,p-q,\mu x)\in W_j \\
0 &\text{otherwise.}
\end{cases}
\]

Then $g$ is continuous, and the support of $g$ is compact since it
is contained in $Z(s(\overline{U_1\cdots U_p})\cap
s(\overline{V_1\cdots V_q}))$. Hence $g\in C_c(G_\Lambda^{(0)})$,
and it follows that $g = \Psi_0\times_\Tt\Psi_1(a)$ for some
$a\in\Tt(E_\Lambda)$.

Define
\[
W = \{(\lambda,\mu)\in (U_1\cdots U_p)*_s(V_1\cdots V_q) : (\lambda
x,p-q,\mu x)\in\supp(\varphi_j f) \text{ for some } x\},
\]
and let $P_1,P_2:\Lambda*_s\Lambda\to\Lambda$ be the coordinate
projections. Then since $W$ is compact and $P_1$ is continuous, for
$i=1,\dots,p$, the set $\Seg^p_{(i-1,i)}(P_1(W))$ is a compact
subset of $U_i$. By Urysohn's Lemma, there exists $\xi_i\in
C_c(U_i)$ satisfying
\[
\xi_i(e) = 1 \quad\text{for all } e\in \Seg^p_{(i-1,i)}(P_1(W)).
\]
Similarly, for $i=1,\dots,q$, there exists $\eta_i\in C_c(V_i)$ such
that
\[
\eta_i(e) = 1 \quad\text{for all } e\in \Seg^q_{(i-1,1)}(P_2(W)).
\]

Calculations using the convolution product on $C_c(G_\Lambda)$
reveal
\[
\varphi_jf =
\Psi_1(\xi_1)\cdots\Psi_1^p(\xi_p)\big(\Psi_0\times_\Tt\Psi_1(a)
\big)\Psi_1(\eta_q)^*\cdots \Psi_1(\eta_1)^*,
\]
so $f$ is in the image of $\Psi_0\times_\Tt \Psi_1$. Since
$\Psi_0\times_\Tt\Psi_1$ has closed range, it follows that
$\Psi_0\times_\Tt\Psi_1$ is surjective.
\end{proof}

\begin{proof}[{Proof of Theorem~\ref{thm:Toeplitz algebra of top
graphs correspond}}]
Fix a faithful nondegenerate representation $\pi$ of
$C^*(G_\Lambda)$ on a Hilbert space $H$. Then
$(\pi\circ\Psi_0,\pi\circ\Psi_1)$ is a Toeplitz $E_\Lambda$-pair,
and the universal property of $\Tt(E_\Lambda)$ gives
\begin{equation}\label{eqn:universal map}
(\pi\circ\Psi_0)\times_\Tt(\pi\circ\Psi_1) = \pi\circ
(\Psi_0\times_\Tt\Psi_1).
\end{equation}
To show $\pi\circ(\Psi_0\times_\Tt\Psi_1)$ is faithful, we use
\cite[Theorem~2.1]{FowR2} which says in our setting that
$(\pi\circ\Psi_0)\times_\Tt(\pi\circ\Psi_1): \Tt(E_\Lambda)\to B(H)$
is faithful if $C_0(\Lambda^0)$ acts faithfully on
$((\pi\circ\Psi_1)(C_d(\Lambda^1))H)^\perp$; that is, if for nonzero
$f\in C_0(\Lambda^0)$, there exists $h\in
((\pi\circ\Psi_1)(C_d(\Lambda^1))H)^\perp$ such that
$(\pi\circ\Psi_0)(f)h\neq 0$.

Fix nonzero $f\in C_0(\Lambda^0)$ and $v\in\Lambda^0$ such that
$f(v)\neq 0$. Let $V\subset\Lambda^0$ be a relatively compact open
neighbourhood of $v$ such that $f(w)\neq 0$ for all $w\in V$, so
$\overline{V}\subset\supp(f)$. Then $U := Z(V)\cap Z(\Lambda^1)^c$
is a relatively compact open neighbourhood of $(v,0,v)$ contained in
$\supp(\Psi_0(f))$; choose $g\in C_c(G_\Lambda^{(0)})$ satisfying
$g(v,0,v)\neq 0$ and $\supp(g)\subset U$. In particular,
$(\Psi_0(f)g)(v,0,v)\neq 0$ and
\begin{equation}\label{eqn:support of g}
\supp(g)\subset \{(w,0,w)\in G_\Lambda^{(0)} : w\in\Lambda^0\}.
\end{equation}

Since $\Psi_0(f)g\neq 0$ and $\pi$ is faithful, there exists $h'\in
H$ such that $\pi(\Psi_0(f)g)h'\neq 0$. Defining $h := \pi(g)h'$ and
using \eqref{eqn:support of g}, we deduce $h\in
((\pi\circ\Psi_1)(C_d(\Lambda^1))H)^\perp$, and we have $(\pi\circ
\Psi_0)(f)h=\pi(\Psi_0(f)g)h'\neq 0$. Thus $C_0(\Lambda^0)$ acts
faithfully on $((\pi\circ\Psi_1)(C_d(\Lambda^1))H)^\perp$, and
\cite[Theorem~2.1]{FowR2} implies
$(\pi\circ\Psi_0)\times_\Tt(\pi\circ\Psi_1)$ is faithful. Therefore,
from \eqref{eqn:universal map}, it follows that
$\Psi_0\times_\Tt\Psi_1$ is injective, which, together with
Proposition~\ref{prop:surjectivity}, completes the proof.
\end{proof}

We now construct a Cuntz-Krieger $E_\Lambda$-pair $\Phi$ on
$C^*(\Gg_\Lambda)$. Since $X_\Lambda\setminus \partial\Lambda$ is an
open invariant subset of $G_\Lambda^{(0)}$, we can regard
$C^*(G_\Lambda|_{X_\Lambda\setminus\partial\Lambda})$ as an ideal in
$C^*(G_\Lambda)$ with quotient $C^*(\Gg_\Lambda)$. Letting
$Q:C^*(G_\Lambda)\to C^*(\Gg_\Lambda)$ be the quotient homomorphism
with kernel $C^*(G_\Lambda|_{X_\Lambda\setminus\partial\Lambda})$,
we define $\Phi_0 := Q\circ\Psi_0$ and $\Phi_1 := Q\circ\Psi_1$.

\begin{prop}
The pair $(\Phi_0,\Phi_1)$ is a Cuntz-Krieger $E_\Lambda$-pair on
$C^*(\Gg_\Lambda)$.
\end{prop}

\begin{proof}
Since $(\Psi_0,\Psi_1)$ is a Toeplitz $E_\Lambda$-pair and $Q$ is a
homomorphism, it follows that $(\Phi_0,\Phi_1)$ is a Toeplitz
$E_\Lambda$-pair. We must show that
\begin{equation}\label{eqn:CK relation}
\Phi_0(f) = \Phi^{(1)}(\phi(f)) \quad\text{for all } f\in
C_0(\Lambda^0_{\rm rg}).
\end{equation}
Since $C_c(\Lambda^0_{\rm rg})$ is dense in $C_0(\Lambda^0_{\rm
rg})$, it suffices to check \eqref{eqn:CK relation} holds for $f\in
C_c(\Lambda^0_{\rm rg})$.

Fix $f\in C_c(\Lambda^0_{\rm rg})$. Then $f\circ r\in
C_c(\Lambda^1)$ and $\phi(f) = \theta(f\circ r)$, where $\theta:
C_b(\Lambda^1)\to\Ll(C_d(\Lambda^1))$ is the injective homomorphism
defined by
\[
(\theta(g)\xi)(e) = g(e)\xi(e).
\]
It then follows from \cite[Lemmas~1.15 and~1.16]{Kat1} that there
exist $l\in\NN$ and $\xi_i,\eta_i\in C_c(\Lambda^1)$ for
$i=1,\dots,l$, such that
\begin{equation}\label{eqn:property 1}
f\circ r = \sum_{i=1}^l\xi_i\overline{\eta_i}, \quad\text{where the
product is defined pointwise},
\end{equation}
\begin{equation}\label{eqn:property 2}
\xi_i(e)\overline{\eta_i(e')} = 0 \quad\text{for all } i \text{ and
} e,e'\in\Lambda^1 \text{ with } s(e)=s(e') \text{ and } e\neq e',
\end{equation}
and
\begin{equation}\label{eqn:property 3}
\phi(f) = \sum_{i=1}^l \xi_i\otimes \eta_i^*.
\end{equation}
Thus we have
\[
\Phi^{(1)}(\phi(f)) = \sum_{i=1}^l\Phi_1(\xi_i)\Phi_1(\eta)^*,
\]
so we must show
\begin{equation}\label{eqn:CK relation modified}
\Phi_0(f)(x,m,y) = \sum_{i=1}^l\Phi_1(\xi_i)\Phi_1(\eta_i)^*(x,m,y)
\end{equation}
for all $(x,m,y)\in \Gg_\Lambda$.

Fixing $(x,m,y)\in \Gg_\Lambda$, the right-hand side of
\eqref{eqn:CK relation modified} is equal to
\begin{equation}\label{eqn:CK relation modified 2}
\sum_{i=1}^l\sum_{(x,n,z)\in\Gg_\Lambda}\Phi_1(\xi_i)(x,n,z)
\overline{\Phi(\eta_i)(y,n-m,z)}.
\end{equation}
A summand from \eqref{eqn:CK relation modified 2} may be nonzero
only if $n=1$, $\sigma^1x=z$, $n-m=1$ and $\sigma^1y=z$; that is,
only if $m=0$, $n=1$, $\sigma^1x=z=\sigma^1y$ and, consequently,
$d(x),d(y)\ge 1$. Thus
\begin{align}
\sum_{i=1}^l&\Phi_1(\xi_i)\Phi_1(\eta_i)^*(x,m,y) \notag\\
&= {\begin{cases}
\delta_{m,0}\delta_{\sigma^1x,\sigma^1y}\sum_{i=1}^l\xi_i(x(0,1))
\overline{\eta_i(y(0,1))} &\text{if } d(x),d(y)\ge 1 \\
0 &\text{otherwise}
\end{cases}} \notag\\
&= {\begin{cases}
\delta_{m,0},\delta_{x,y}\sum_{i=1}^l\xi_i(x(0,1))
\overline{\eta_i(x(0,1))} &\text{if } d(x)\ge 1 \\
0 &\text{otherwise}
\end{cases}}\quad\text{by \eqref{eqn:property 2}}\notag \\
&= {\begin{cases}
\delta_{m,0}\delta_{x,y}(f\circ r)(x(0,1)) &\text{if } d(x)\ge 1 \\
0 &\text{otherwise}
\end{cases}} \quad\text{by \eqref{eqn:property 1}}\notag \\
&= {\begin{cases}
\delta_{m,0}\delta_{x,y}\Phi_0(f)(x,0,x) &\text{if } d(x)\ge 1 \\
0 &\text{otherwise}
\end{cases}}\label{eqn:almost there}
\end{align}
Hence, if $d(x)\ge 1$ then \eqref{eqn:CK relation modified} follows
from \eqref{eqn:almost there}. On the other hand, if $d(x)=0$ then
$x=s(x)\not\in \Lambda^0_{\rm rg}$, and since
$\supp(f)\subset\Lambda^0_{\rm rg}$, it follows that
$\Phi_0(f)(x,0,x) = 0$. Therefore \eqref{eqn:CK relation modified}
holds for all $(x,m,y)\in\Gg_\Lambda$, and $(\Phi_0,\Phi_1)$ is a
Cuntz-Krieger $E_\Lambda$-pair.
\end{proof}

\begin{proof}[{Proof of Theorem~\ref{thm:C-K algebras of top graphs
correspond}}]
Since $\Phi$ is a Cuntz-Krieger $E_\Lambda$-pair, it is also a
Toeplitz $E_\Lambda$-pair, and we have two induced homomorphisms
\[
\Phi_0\times_{\Tt}\Phi_1 :\Tt(E_\Lambda)\to C^*(\Gg_\Lambda)
\quad\text{and}\quad \Phi_0\times_{\Oo}\Phi_1 : \Oo(E_\Lambda)\to
C^*(\Gg_\Lambda).
\]

By considering the universal Cuntz-Krieger $E_\Lambda$-pair
$(j_0,j_1)$ on $\Oo(E_\Lambda)$, the universal property of
$\Tt(E_\Lambda)$ gives a homomorphism
$j_0\times_{\Tt}j_1:\Tt(E_\Lambda)\to\Oo(E_\Lambda)$ which
satisfies
\[
\Phi_0\times_{\Tt}\Phi_1 = (\Phi_0\times_{\Oo}\Phi_1)\circ
(j_0\times_{\Tt}j_1).
\]
On the other hand, recalling the quotient homomorphism
$Q:C^*(G_\Lambda)\to C^*(\Gg_\Lambda)$, the universal property of
$\Tt(E_\Lambda)$ gives
\[
\Phi_0\times_{\Tt}\Phi_1 = Q\circ (\Psi_0\times_{\Tt}\Psi_1).
\]
Hence
\[
(\Phi_0\times_{\Oo}\Phi_1)\circ (j_0\times_{\Tt}j_1) = Q\circ
(\Psi_0\times_{\Tt}\Phi_1),
\]
and surjectivity of $\Phi_0\times_{\Oo}\Phi_1$ follows from
surjectivity of $Q$ and $\Psi_0\times_{\Tt}\Psi_1$.

To show that $\Phi_0\times_\Oo\Phi_1$ is injective, we use the
gauge-invariant uniqueness theorem \cite[Proposition~4.5]{Kat1}.

The map $c:\Gg_\Lambda\to\ZZ$ defined by $c(x,m,y) = m$ is a
continuous functor, hence by \cite[Proposition~II.5.1]{Ren1} there
is a strongly continuous action $\beta$ of $\widehat{\ZZ}=\TT$ on
$C^*(\Gg_\Lambda)$ such that $\beta_t(g)(x,m,y) = t^mg(x,m,y)$ for
all $t\in\TT$ and $g\in C_c(\Gg_\Lambda)$, and such that $\beta$
leaves $C_0(\Gg_\Lambda^{(0)})$ pointwise fixed. We then have
$\beta_t(\Phi_0(f))=\Phi_0(f)$ and $\beta_t(\Phi_1(\xi)) =
t\Phi_1(\xi)$ for all $t\in\TT$, $f\in C_0(\Lambda^0)$ and $\xi\in
C_d(\Lambda^1)$.

It is straightforward to see that $\Phi_0:C_0(\Lambda^0)\to
C^*(\Gg_\Lambda)$ is injective. Therefore
\cite[Proposition~4.5]{Kat1} implies
$\Phi_0\times_\Oo\Phi_1:\Oo(E_\Lambda)\to C^*(\Gg_\Lambda)$ is
injective, and we're done.
\end{proof}

\end{document}